\documentclass[10pt]{article}
\usepackage[utf8]{inputenc}
\usepackage[english]{babel}
\usepackage{amsmath,amsthm,amsfonts}
\usepackage{graphicx}
\usepackage{bm}
\usepackage{listings}
\usepackage{multirow}
\usepackage{color}

\graphicspath{ {./figs/} }

\newcommand{\grad}{\mathop{\rm grad}\nolimits}
\renewcommand{\div}{\mathop{\rm div}\nolimits}
\newcommand{\const}{\mathop{\rm const}\nolimits}

\title{Machine learning for accelerating effective property prediction for poroelasticity problem in stochastic media}

\author{
Maria Vasilyeva \thanks{Institute for Scientific Computation, Texas A\&M University, College Station, TX 77843-3368 \& Multiscale model reduction laboratory, North-Eastern Federal University, Yakutsk, Republic of Sakha (Yakutia), Russia, 677980. Email: {\tt vasilyevadotmdotv@gmail.com}.}
\and
Aleksey Tyrylgin 
\thanks{Multiscale model reduction laboratory, North-Eastern Federal University, Yakutsk, Republic of Sakha (Yakutia), Russia, 677980.}
}

\begin{document}

\maketitle

\begin{abstract}
In this paper, we consider a numerical homogenization of the poroelasticity problem with stochastic properties.
The proposed method based on the construction of the deep neural network (DNN) for fast calculation of the effective properties for a coarse grid approximation of the problem.
We train neural networks on the set of the selected realizations of the local microscale stochastic fields and macroscale characteristics (permeability and elasticity tensors).
We construct a deep learning method through convolutional neural network (CNN) to learn a map between stochastic fields and effective properties. Numerical results are presented for two and three-dimensional model problems and show that proposed method provide fast and accurate effective property predictions.
\end{abstract}

\section{Introduction}
% stochastic media, uncertaintly
Uncertainties remain in most models of real world problems and arise due to lack of knowledge of heterogeneous properties.
Uncertainties may be described by stochastic models with uncertain parameters . %that described by random variables.
% subsurface
In the subsurface processes an uncertainty of the properties varying in space. Numerical solution of such problems is difficult and some type of coarsening is necessary to perform fast calculations without resolving small scale heterogeneity by grid.
However, most of coarse grid methods are developed for fixed realization of the property field and such coarse grid models may not be sufficient for fast simulations.

% coarse grid
Classical numerical methods for solving problems with heterogeneous properties are based on the standard finite element or finite volume approximations.
For the applicability and convergence, the size of the computational grid must be small enough to explicitly resolve all the heterogeneity by the grid.
Approximation on a fine grid significantly increases the dimension of the discrete problem and requires large computational resources.
For such problems, we should use multiscale or upscaling methods to create effective approximations on a coarse grid \cite{vashomog, vasilyeva2016tp, brown2016generalized, brown2016generalized2, vasilyeva2018nlmcporoel, akkutlu2018gmsfemporoel}.
% homog
For the periodic media, an asymptotic two-scale homogenization method can be used, where coarse grid (macroscale) equations with effective medium properties are derived. In this method an effective property are calculated for one heterogeneity period by solving local problems \cite{allaire1992homogenization, sanchez1980non, bakhvalov1984homogenization}.
For problems in non-periodic media, the methods of numerical homogenization (upscaling) are used, where local problems are solved to calculate effective characteristics in each local domains.
For highly accurate solution a various multiscale methods are developed. Multiscale and homogenization methods assume the solution of local problems for taking into account microscale information for approximation on a macroscale computational grid. Solving such problems are computationally expensive and require calculations for each realization of the property field.
Therefore, accelerating of the local calculations is necessary for fast simulations \cite{pilania2013accelerating, aarnes2008mixed, ganis2008stochastic, dostert2006coarse, efendiev2012systematic, he2013stochastic}.

% ml
In recent years, many new highly effective methods for constructing machine learning algorithms are developed. The reason of the increased usage and popularity in recent times associated with development of easy-to-use open source software libraries and availability of graphics processing units (GPUs) for accelerated computation. Furthermore, this jump is associated with development of the deep learning, advance stochastic optimization techniques and robust regularization techniques such as dropout \cite{kingma2014adam, srivastava2014dropout, lecun1998gradient, lecun2015deep}.
Deep neural network dramatically improved accuracy of the machine learning methods due to presence of the multiple processing layers that learn representations of data with multiple levels of abstraction (feature extraction) \cite{lecun2015deep}. With the composition of enough layers, very complex functions can be learned.
Convolutional neural network is a one particular type of deep network that was much easier to train and generalized much better than networks with full connectivity between layers \cite{krizhevsky2012imagenet, lecun2015deep, simonyan2014very}.
The architecture of a convolutional neural network contains composing of the convolutional and pooling layers that passed through a non-linear layers. Several layers of convolution, non-linearity and pooling are stacked, and completes with a several fully-connected layers. The key aspect of deep learning is that these layers of features are not designed by human engineers: they are learned from data using a general-purpose learning procedure.
Nowadays, convolutional neural networks have become very popular and are used to solve various problems, including calculation of the physical properties, identifying new dependencies, and solution  of the partial differential equations without direct calculations of grid problems \cite{srisutthiyakorn2016deep, sudakov2018driving, chan2018machine}.

% ml in this work
In this work, we construct a machine learning algorithm for an accurate and fast calculations of the effective properties of the random poroelastic media (permeability and elasticity tensors).
We develop a machine learning method through a construction of convolutional neural network (CNN) to learn a map between  stochastic fields and effective properties.
For training convolutional network, we calculate reference effective properties by the solution of the local problems, whereas input data we use a local heterogeneous property on the structured grid (array of pixel values).
After that, we use this network to significantly reduce the calculation time of effective characteristics and coarse grid solution of the poroelasticity problem in stochastic media.

The paper is organized as follows. In Section 2, we consider model poroelasticity problem in stochastic media with fine grid approximation. Next, we present a numerical homogenization technique for coarse grid solution in Section 3. In Section 4, we present a machine learning algorithm for accelerating effective property prediction for stochastic poroelastic media. In Section 5, we present numerical results for two and three-dimensional problems for random media with exponential variogram, where we present errors of numerical homogenization for some test cases, results of the training of the machine learning algorithms,  relative errors of algorithm for several samples of realizations. Finally, in Section 5,  we present and discuss time of machine learning algorithm construction and prediction time.

\section{Model problem}

Let $p$ and $u$ are the pressure and displacement. In domain $\Omega$, we  consider linear poroelasticity problem for $(p, u)$ \cite{vassplitting, coussy2004poromechanics, kim2011stability2}
\begin{equation}
\label{eq:1}
\begin{split}
\frac{1}{M} \frac{\partial p}{\partial t} +  \alpha \frac{\partial \div  u}{\partial t}  - \div \left(\frac{k(x, \omega)}{\nu_f}  \grad p \right) &= f, \quad  x \in \Omega,\\
- \div (\sigma(u) ) + \alpha \grad p  &= 0, \quad  x  \in \Omega, 
\end{split}
\end{equation}
where $C$ is  the elasticity tensor, $k$ is the permeability, $\sigma$ is the stress tensor, $\nu_f$ is the fluid viscosity, $f$ is the source term, $M$ is the Biot modulus and $\alpha$ is the Biot-Willis fluid-solid coupling coefficient and  
\[
 \sigma(u) = C(x, \omega) : \varepsilon(u),
\quad
 \varepsilon(u) = \frac{1}{2} \left( \grad  u + \grad  u^T \right).
\]
We note that, $k$ and $C$ are stochastic coefficients  and  $M$, $\alpha$, $\nu_f$  are constants. 
Since the permeability and elastic modulus are a stochastic function, $p$ and $u$ are also stochastic. 

Here for two-dimensional case (2D), we have 
\[
u = (u_1, u_2), \quad 
\sigma = (\sigma_1, \sigma_2, \sigma_{12})^T, \quad  
\varepsilon = (\varepsilon_1, \varepsilon_2, \varepsilon_{12})^T,
\] 
\begin{equation}
k =
\begin{bmatrix}
k_{11} & k_{12}\\
k_{21} & k_{22}
\end{bmatrix}, 
\quad 
C = \begin{bmatrix}
C_{1111}  & C_{1122} &  C_{1112} \\
C_{2211}  & C_{2222} &  C_{2212} \\
C_{1211}  & C_{1222} &  C_{1212}
\end{bmatrix},
\end{equation}
and for three-dimensional case (3D)
\[
u = (u_1, u_2, u_3), \quad 
\sigma = (\sigma_1, \sigma_2, \sigma_3,
 \sigma_{12}, \sigma_{23}, \sigma_{13})^T, \quad  
\varepsilon = (\varepsilon_1, \varepsilon_2, \varepsilon_3, 
\varepsilon_{12}, \varepsilon_{23}, \varepsilon_{13})^T,
\] 
\begin{equation}
k =
\begin{bmatrix}
k_{11} & k_{12} & k_{13}\\
k_{21} & k_{22} & k_{23}\\
k_{31} & k_{32} & k_{33}\\
\end{bmatrix}, 
\quad
C = \begin{bmatrix}
C_{1111}  & C_{1122} &  C_{1133} &  C_{1112} &  C_{1123} &  C_{1131} \\
C_{2211}  & C_{2222} &  C_{2233} &  C_{2212} &  C_{2223} &  C_{2231} \\
C_{3311}  & C_{3322} &  C_{3333} &  C_{3312} &  C_{3323} &  C_{3331} \\
C_{1211}  & C_{1222} &  C_{1233} &  C_{1212} &  C_{1223} &  C_{1231} \\
C_{2311}  & C_{2322} &  C_{2333} &  C_{2312} &  C_{2323} &  C_{2331} \\
C_{3111}  & C_{3122} &  C_{3133} &  C_{3112} &  C_{3123} &  C_{3131} \\
\end{bmatrix}. 
\end{equation}

%bc
We consider system of equations \eqref{eq:1} with following  initial and boundary conditions 
\[
p = p_0, \quad x \in \Omega, \quad t = 0, 
\]\[
p = p_1, \quad x \in \Gamma_p, \quad \text{ and } \quad
\frac{\partial p}{\partial n} = 0, \quad x \in \partial \Omega / \Gamma_p,
\]\[
u = 0, \quad  \quad x \in \Gamma_u, \quad \text{ and } \quad 
\sigma \cdot n = 0, \quad x \in \partial \Omega / \Gamma_u. 
\]

\textbf{Variational formulation. }
We use a finite element method to find an approximate solution of the poroelasticity problem. 
Let
\[
V = \{ v \in H^1(\Omega): v = p_1, \, x \in  \Gamma_p \}, \quad 
\hat{V} = \{ v \in H^1(\Omega): v = 0, \, x \in  \Gamma_p \}.
\] \[
W = \{ w \in [H^1(\Omega)]^d: w = 0, \, x \in  \Gamma_u \}, \quad 
\hat{W} = W.
\] 
where $d = 2,3$. 
We have following variational formulation of the problem: find  $(u, p) \in V \times Q$ such that
\begin{equation}
\label{eq:fem}
\begin{split}
d \left( \frac{\partial  u}{\partial t}, q \right)  + 
m \left( \frac{\partial p}{\partial t}, q \right) + 
b(p , q) &= l(q), \quad \forall  q \in  \hat{Q}, \\
a( u,  v) + g(v, p) &= 0, \quad \forall  v \in \hat{V}, 
\end{split}
\end{equation}
where
\[
a(u,  v)  = \int_{\Omega}  \sigma(u) : \varepsilon(v) \, dx, 
\quad
b(p, q) = \int_{\Omega} \frac{k}{\nu_f} \grad p \, \cdot \, \grad q \, dx,
\quad
m(p, q) = \int_{\Omega} \frac{1}{M} \, p \, q \, dx,  
\]\[
d(u, p) = \int_{\Omega} \alpha \, \div  u \, p \, dx , 
\quad
g(v, p) = \int_{\Omega} \alpha \, v \cdot \grad p \, dx , 
\quad
l(q) = \int_{\Omega} f  \, q \, dx.
\]

\textbf{Fine grid system. }
Let $\mathcal{T}^h$ is the fine grid partition of the computational domain $\Omega$ into finite elements.  In particular, we use piecewise linear basis functions for finite element approximation. The standard implicit finite difference scheme is used for the approximation with time step size $\tau$ and superscripts $n$, $n+1$ denote previous and current time levels. 
We will use fine grid formulation for reference  solution and error calculations in Section 5. 
On the fine grid, the equation (\ref{eq:fem}) can be presented in matrix form:
\begin{equation}
\label{eq:mat}
\begin{split}
D \frac{u^{n+1} - u^n}{\tau} + M \frac{p^{n+1} - p^n}{\tau} 
+ B p^{n+1} & = F,\\
A u^{n+1} + G p^{n+1}& = 0,
\end{split}
\end{equation}
where  
$M = [m_{ij}], \quad m_{ij} =  m(\phi_i, \phi_j)$,  
$B = [b_{ij}], \quad b_{ij} =  b(\phi_i, \phi_j)$,  
$A = [a_{ij}], \quad a_{ij} =  a(\psi_i, \psi_j)$,  
$D = [d_{ij}], \quad d_{ij} =  d(\psi_i, \phi_j)$, 
$G = [g_{ij}], \quad g_{ij} =  g(\psi_i, \phi_j)$     
and $F = \{ f_j \}$, $f_j  = l(\phi_j)$,  $\phi_i$ and $\psi_i$ are the linear  basis functions defined on $\mathcal{T}^h$. % $i,j = 1...N_f$ and $N_f$ denotes the number of the nodes on the fine grid.

\section{Numerical homogenization}

The size of the computational grids for heterogeneous media must be small enough to explicitly resolve all the heterogeneous by the grid.  
Approximation on a fine grid significantly increases the dimension of the discrete problem, therefore the computational time, as well as the amount of used memory are also increase. 
For construction of the coarse grid approximation, we use a numerical homogenization technique, 
where we solve local problems on a fine grid  to identify effective coefficients for each coarse grid cell \cite{vashomog, cottereau2013numerical, hazanov1994order, moravec2009numerical}. 

Let $\mathcal{T}^H = \cup_i K_i$ ($i = \overline{1,N_c}$) be a structured partition of the computational domain $\Omega$ into elements $K_i$, where $N_c$ is the number of the coarse grid cells and $i$ is the coarse grid cell index.
We calculate effective permeability and elastic coefficient, $k^*$ and $C^*$ in each coarse cell for nonperiodic heterogeneous media for some given realization.%, $\{C^{\star, K_i}(\omega)\}_{i=1..N}$ and $\{k^{\star, K_i}(\omega)\}_{i=1..N}$ for $K_i \in \mathcal{T}_H$.

\textbf{Permeability tensor. }
For calculation of the effective permeability, we solve following local problem in $K_i$
\begin{equation}
\label{eq:k1}
\begin{split}
- \div \left( k^{K_i}(x, \omega) \grad \psi^{K_i}_j \right) &= 0,
\quad x \text{ in } K_i, \\
\psi^{K_i}_j &= x_j , 
\quad x \text{ on } \partial K_i, 
\end{split}
\end{equation}
where $k^{K_i} $ is the restriction of the heterogeneous coefficient $k(x, \omega)$ to local domain $K_i$. 
Here $x = (x_1, x_2)$ for 2D case and $x = (x_1, x_2, x_3)$ for 3D case. Therefore, for 2D problem, we solve two local problems and for 3D problem, we have three local problems.
Note that, another boundary conditions can be applied for local problem.

Next, we can find elements of the effective permeability tensor for current $K_i$
\begin{equation}
\label{k2}
k^{\star, K_i}_{lj}(\omega) = 
\frac{1}{|K_i|}
\int_{K_i}  k^{K_i}(x, \omega) \frac{\partial \psi^{K_i}_l}{\partial x_j}  dx, 
\quad l,j = \overline{1,d},
\end{equation}
where $d$ is problem dimention, $d = 2$ or $3$.  This permeability tensor is symmetric. 

\textbf{Elasticity tensor. }
For effective elastic modulus, we apply similar algorithm and solve following local problem in $K_i$
\begin{equation}
\label{E1}
\begin{split}
- \div  (C^{K_i}(x, \omega) : \varepsilon(\phi_{rs}^{K_i}) ) = 0, \quad x \text{ in } K_i,  \\
\phi_{rs}^{K_i} = \Lambda^{(rs)} x, \quad x \text{ on } \partial K_i,
\end{split}
\end{equation}
where $\phi^{K_i} = (\phi_1^{K_i}, \phi_2^{K_i})$ for $d = 2$, $\phi^{K_i} = (\phi_1^{K_i}, \phi_2^{K_i}, \phi_3^{K_i})$ for $d = 3$ and 
\[
 \Lambda_{ij}^{(rs)} = \frac{1}{2}
 \left( \delta_{ir} \delta_{js} + \delta_{is} \delta_{jr} \right),  \quad r, s = \overline{1,d}.
\]

The elements of the effective elastic modulus are calculated as follows
\begin{equation}
\label{E2}
C^{\star, K_i}_{rspq}(\omega) = \frac{1}{|K_i|}
\int_{K_i} C(x, \omega) \varepsilon (\phi^{(rs)}) : \varepsilon (\phi^{(pq)})\, dx, \quad r,s,p,q = \overline{1,d}.
\end{equation}
This elaticity tensor is symmetric. 

There are existed several approaches for the numerical homogenization methods based on the two-scale asymptotic analysis with solution of the local problems in representative volume $K_i$ with periodic boundary condition, Dirichlet or Neumann boundary condition, or mixed boundary condition.

\textbf{Coarse grid system. }
Finally, we solve the poroelasticity problem on a coarse grid with precalculated effective permeability and elastic modulus
\begin{equation}
\label{eq:2}
\begin{split}
\frac{1}{M} \frac{\partial p}{\partial t} +  \alpha \frac{\partial \div  u}{\partial t}  - \div \left( \frac{k^{\star}(x, \omega)}{\nu_f} \grad p \right) &= f, \\
- \div \left( C^{\star}(x, \omega) : \varepsilon (u) \right) + \alpha \grad p  &= 0, 
\end{split}
\end{equation}
using Galerkin finite element method. 
We note that, for each random field realization, we should calculate effective properties. This calculation can be computationally expensive due to solution of the local problems in each coarse cell up to fine grid resolution of the heterogeneous fields. 
Next, we will describe the construction of the machine learning algorithm for the fast calculation. We will use presented properties calculation technique for creating a dataset for train and test a deep neural network.

\section{Machine learning algorithm}

In this section, we  present a machine learning approach for solution of the stochastic poroelasticity problem.  The main idea is to use different realizations of the permeability and elastic coefficient fields to train and test a deep neural network. After that, constructed deep network is used to fast and accurate solution of the coarse grid poroelasticity system. 

We have following main steps:
\begin{itemize}
\item Generate dataset to train and test of the neural network. 
\item Construction of the neural network and test it on a set of realizations.
\item Fast construction and solution of the coarse grid system using trained neural network for effective property prediction. 
\end{itemize}

Next, we consider generation of the dataset and network construction  in details. 

\subsection{Stochastic properties}

We suppose  isotropic permeability field and  for elasticity tensor, we have
\[
C = \begin{bmatrix}
\lambda + 2 \mu & \lambda & 0 \\
\lambda & \lambda + 2 \mu & 0 \\
0 & 0 & 2 \mu
\end{bmatrix},
\quad 
C = \begin{bmatrix}
\lambda + 2 \mu & \lambda & \lambda & 0 & 0 & 0 \\
\lambda & \lambda + 2 \mu & \lambda & 0 & 0 & 0\\
\lambda & \lambda & \lambda + 2 \mu & 0 & 0 & 0\\
0 & 0 & 0 & 2 \mu & 0 & 0 \\
0 & 0 & 0 & 0 & 2 \mu & 0 \\
0 & 0 & 0 & 0 & 0 & 2 \mu \\
\end{bmatrix},
\]
for 2D and 3D problems, respectively. 
Here for Lame parameters, we have 
\[
\mu = \frac{E}{2 (1 + \eta)}, \quad
\lambda = \frac{E \eta}{(1+ \eta) ( 1- 2 \eta)}.
\]
where $\eta = \const$ is the Poisson's ratio.

Let $Y^{}(x,\omega)$ is the random normal field  generated by the Karhunen-Lo{\'e}ve (KL) expansion with the corresponding covariance matrix $R(x,y)$  \cite{aarnes2008mixed, ganis2008stochastic}.  
We suppose that the Covariance structure $R(x,y)$ is of the form
\begin{align}
R(x,y)=\sigma^2\exp\left(- \sum_{i=1}^{d} \frac{|x_i - y_i|^2}{2 l_i^2} \right), 
\end{align}
where  $l_1$, $l_2$, $l_3$ are the correlation lengths and $\sigma^2$ is the variance.

We have
\[
 Y(x,\omega)=\sum_{k=1}^{N_Y}\sqrt{\lambda_k}\nu_k(\omega)\varphi_k(x),
\]
where $\nu_k$ is the stochastic coefficients, $\varphi_k$ and  $\lambda_k$  are the eigenfunctions and eigenvalues . 

The stochastic permeability and elastic modulus fields  ($\kappa = \kappa(x, \omega)$, $E = E(x, \omega)$) are given as follows
\begin{equation}
\label{kandE}
k(x,\omega)=\exp\left(Y^{}(x,\omega)\right), \quad 
E(x,\omega)=\bar{E}(x)+\alpha Y(x,\omega),
\end{equation}
where $\alpha>0$ is the strength of the randomness.

\subsection{Dataset}

We construct deep neural network, whereas input parameters we use a set of fine grid random properties in local domains and as output, we use an effective properties ($k^{\star}$ and $C^{\star}$).

We have following steps:
\begin{itemize}
\item Generate a set of random realizations of ${Y(x, \omega_l)}$ ($l = \overline{1,M}$) and  
\[
k(x,\omega_l)=\exp\left(Y^{}(x,\omega_l)\right),
\quad 
E(x,\omega_l)=\bar{E}(x)+\alpha Y(x,\omega_l).
\]
\item Define a uniform $N \times N$ coarse grid $\mathcal{T}^H = \cup_i K_i$ with $i = \overline{1,N_c}$ and $N_c = N \cdot N$.
\item For each realization $k(x,\omega_l)$:
\begin{itemize}
\item Divide permeability into local domains, $k^{K_i}(x,\omega_l)$.
\item Calculate effective permeability tensor $k^{*,K_i}(\omega_l)$ by solution of the local problems in $K_i$. 
\end{itemize}
\item For each realization $E(x,\omega_l)$:
\begin{itemize}
\item Divide elastic modulus into local domains, $E^{K_i}(x,\omega_l)$.
\item Calculate effective elasticity tensor $C^{*,K_i}$ by solution of the local problems in $K_i$. 
\end{itemize}
\item Save $\{k^{K_i}(x,\omega_l) \rightarrow  k^{*,K_i}(\omega_l)\}$ and  $\{E^{K_i}(x,\omega_l) \rightarrow  C^{*,K_i}(\omega_l)\}$ for $i = \overline{1,N_c}$ and $l = \overline{1,M}$.
\item Normalize dataset and use it for machine learning algorithm construction. 
\end{itemize}

We have two datasets: (1)  2D stochastic fields and (2) 3D stochastic fields. 
Each of the stochastic field is represented as $d$-dimensional array with size $N_l^d$ ($d = 2,3$), where each value represent normalized property of the heterogeneous media. 
The scale of each array in dataset  are scaled to fall within the range $0$ to $1$. 
For representing heterogeneous properties, the high resolution representation should be used for better accuracy of the homogenization methods. 

The input of the network is  $d$ - dimensional normalized array 
\[
X_j = k^{K_i}(x,\omega_l) \quad  \text{ or } \quad 
X_j = E^{K_i}(x,\omega_l), \quad 
j = l \cdot M + i, \quad 
i = \overline{1,N_c}, \quad 
l = \overline{1,M}, 
\]
where $X_j$ has size $N_l^d$ and $j = \overline{1,L}$ is the size of the dataset, where $L = M \cdot N_c$. 

The output for the numerical homogenization is the  normalized array of the effective properties%, where each property normalized separately
\[
Y_j = k^{*,K_i}(\omega_l) \quad  \text{ or } \quad 
Y_j = C^{*,K_i} (\omega_l), \quad 
j = l \cdot M + i, \quad 
i = \overline{1,N_c}, \quad 
l = \overline{1,M}, 
\]
where $k^{*,K_i}$ and $C^{*,K_i}$ are tensors that represented as array with symmetric property. 

The dataset  $\{X_j \rightarrow Y_j\}$ with $j = \overline{1,L}$ contains $L = N_c \cdot M$ local stochastic fields and  used to train, validate and test the neural network. 
In general, effective properties (output) can be obtained from laboratory measurement and input dataset can be obtained using high resolution tomography \cite{srisutthiyakorn2016deep, sudakov2018driving, chan2018machine}. For example for permeability calculation, digital rock database can be used. % and single phase fluid flow with Navier-Stokes equations is the established method for solving absolute permeability. 

\subsection{Network}

In recent years, many new highly effective methods for constructing artificial intelligence have appeared. 
This jump is associated with the development of new methods based on the construction of deep learning methods through convolutional neural networks, and the development of open source libraries that accelerate calculations on GPUs.  
The architecture of a convolutional neural network (CNN) contains composing of the convolutional and pooling layers that passed through a non-linearity such as a ReLU functions \cite{lecun2015deep, krizhevsky2012imagenet}.  Several layers of convolution, non-linearity and pooling are stacked, and completes with a several fully-connected layers. 
%This network was first proposed in a widely known work in which the foundations of the architectures of modern convolutional neural networks are presented.

We train a convolutional neural network by a dataset $\{X_j \rightarrow Y_j\}$ of local random coefficients ($X_j$) and macroscale characteristics (effective medium properties, $Y_j$). 
Constructed machine learning algorithm efficiently determine of dependencies and used  for fast calculation of a effective properties of random media. 

To generate a dataset for machine learning, the samples $X_j$ are generated randomly and to compute train set $Y_j$, we use numerical homogenization 
\[
\{X_1, ..., X_i, ... , X_N \} 
\rightarrow
\text{ Numerical homogenization }
\rightarrow
\{Y_1, ..., Y_i, ... , Y_N \}. 
\]

Next, we divide dataset into train, validation and test sets with sizes $N_{train}$, $N_{val}$ and $N_{test}$. As test set, we take 60 \% of data, another 40 \% randomly divided between train and validation set in 80/20 proportion. 
Therefore, following dataset  is used for train convolutional neural network (CNN)
\[
\text{ Train set: }
\{ (X_1, Y_1), ..., (X_i, Y_i), ... , (X_{N_{train}}, Y_{N_{train}}) \} 
\rightarrow
\, \text{CNN} \,
\]
We use validation set to validate a training process and test data for testing constructed machine learning algorithm. 
After that we use the constructed network as a “black box” for fast prediction of the effective properties for a given local heterogeneity:
\[
X \rightarrow
\, \text{CNN} \,
\rightarrow  Y.
\]

Convolutional neural network is the  deep neural network with  multiple layers that compute a function $F(X_i, W)$, where $X_i$ is the input data and $W$ is the system parameters. By a training process, a machine learning algorithm solve the optimization problem to find model weights that best describe  the train set by minimization of the loss function. 
In this work, we use the mean square error (MSE) as a loss function 
\[
E_{train} = 
\frac{1}{N_{train}} \sum_{i=1}^{N_{train}} |Y_i - F(X_i, W) |^2
\]
and also calculate validation set loss function. For solution of the minimization problem, we use gradient-based optimizer Adam \cite{kingma2014adam}. 

The convolutional neural network has several convolutional and pooling layers with rectified linear units (RELU) activation layer, and two fully connected layers with Dropout.
Dropout algorithm consists of randomly dropping out units of the network during the optimization iteration.
Convolution layers extract and combine local features in $d$ - dimensional arrays due to convolution process that alternates with pooling layers for reduction of the spatial resolution. The pooling layer perform a local averaging. Several layers of convolutions and pooling are alternated in order to detect higher order features for better accuracy of the method.
Successive layers of convolutions and pooling are alternated. Finally, several fully connected layers are applied.

\section{Results}

The numerical calculations of the effective properties has been implemented with the open-source finite element software FEniCS \cite{fenics}. 
Implementation of the machine learning method is based on the open source library Keras \cite{keras}. 
Keras is a open-source library that provide high-level building blocks for developing deep-learning models and based on the several backend engines. We use TensorFlow backend \cite{tensorflow}.

We consider three test cases with different KL-expansion parameters for dataset generation:
\begin{itemize}
\item[] \textit{Test 1}. 2D problem with 
$l_1^2 = 0.2$, $l_2^2 = 0.2$ and $\sigma^2 = 2$. 
\item[] \textit{Test 2}. 2D problem with 
$l_1^2 = 0.1$, $l_2^2 = 0.4$ and $\sigma^2 = 2$. 
\item[] \textit{Test 3}. 3D problem with 
$l_1^2 = 0.2$, $l_2^2 = 0.2$, $l_3^2 = 0.2$ and $\sigma^2 = 2$. 
\end{itemize}

At first, we present results for numerical homogenization method for poroelasticity problem in the stochastic media and compare errors between coarse grid solution and reference (fine grid) solutions. 
Next, we present results for the machine learning algorithm and calculate errors for train and test datasets. 
Finally, we consider a coarse grid solution of the problem, where effective properties are calculated using constructed machine learning method and discuss computational efficiency of the presented method. 

\subsection{Numerical homogenization results}

In this section, we present numerical results of the numerical homogenization method for poroelasticity problems in heterogeneous media. 
For poroelastic problem, we set $\alpha = 1$ and  $M = 1$. 
As initial conditions, we set $p_0 = 0$ and perform calculations for $T_{max} = 0.001$ with 20 time steps. 
We perform calculations using structured coarse and fine grids. 
For each test problems, we consider three cases by varying a heterogeneous elasticity modulus and permeability.

\begin{figure}[h!]
\centering
Test 1\\
\includegraphics[width=0.7\linewidth]{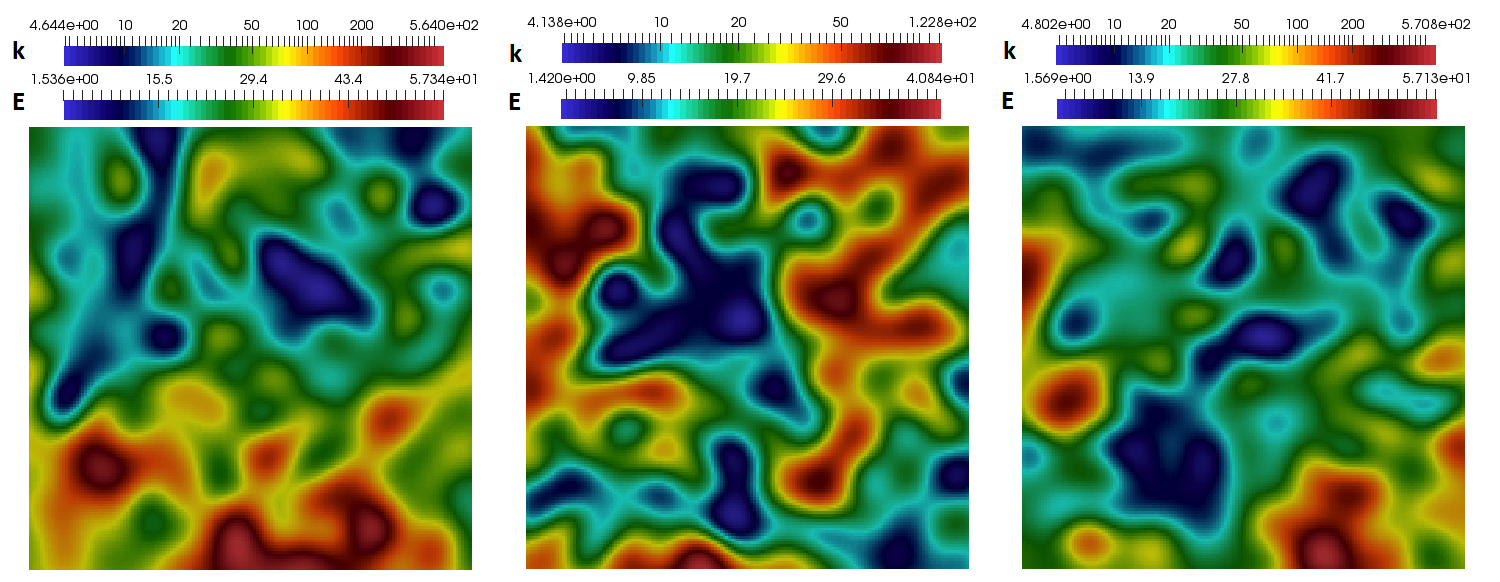} \\
Test 2\\
\includegraphics[width=0.7\linewidth]{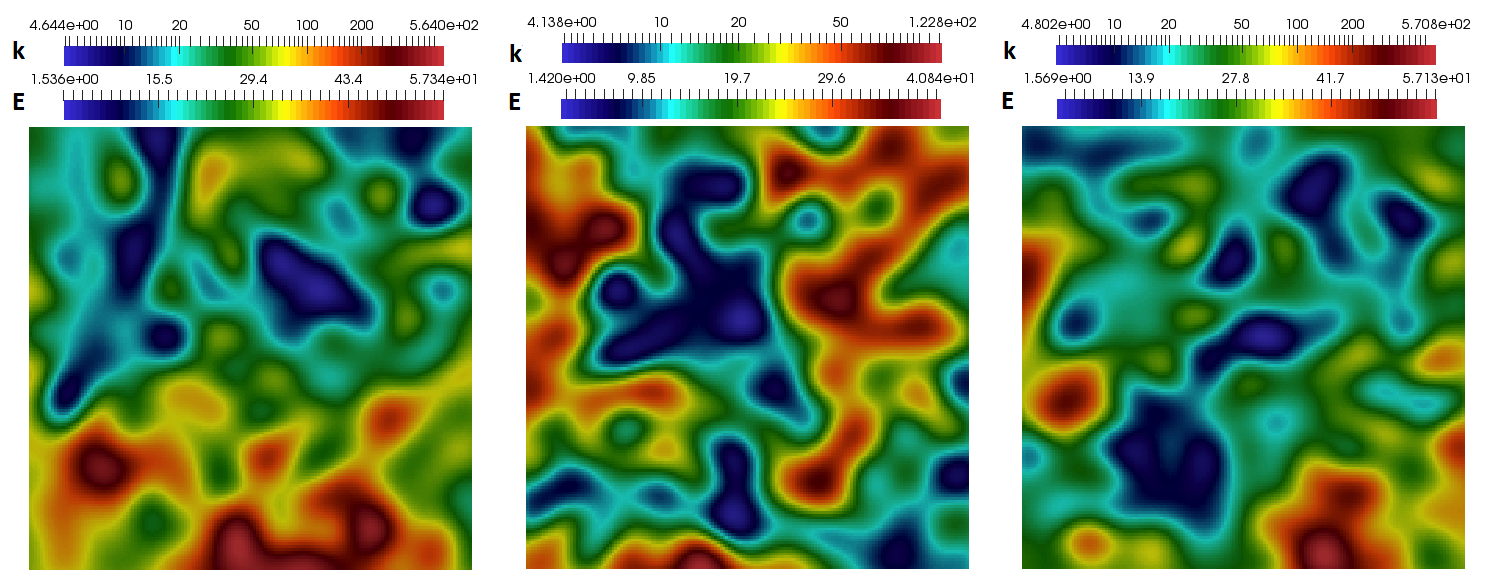} \\
Test 3\\
\includegraphics[width=0.7\linewidth]{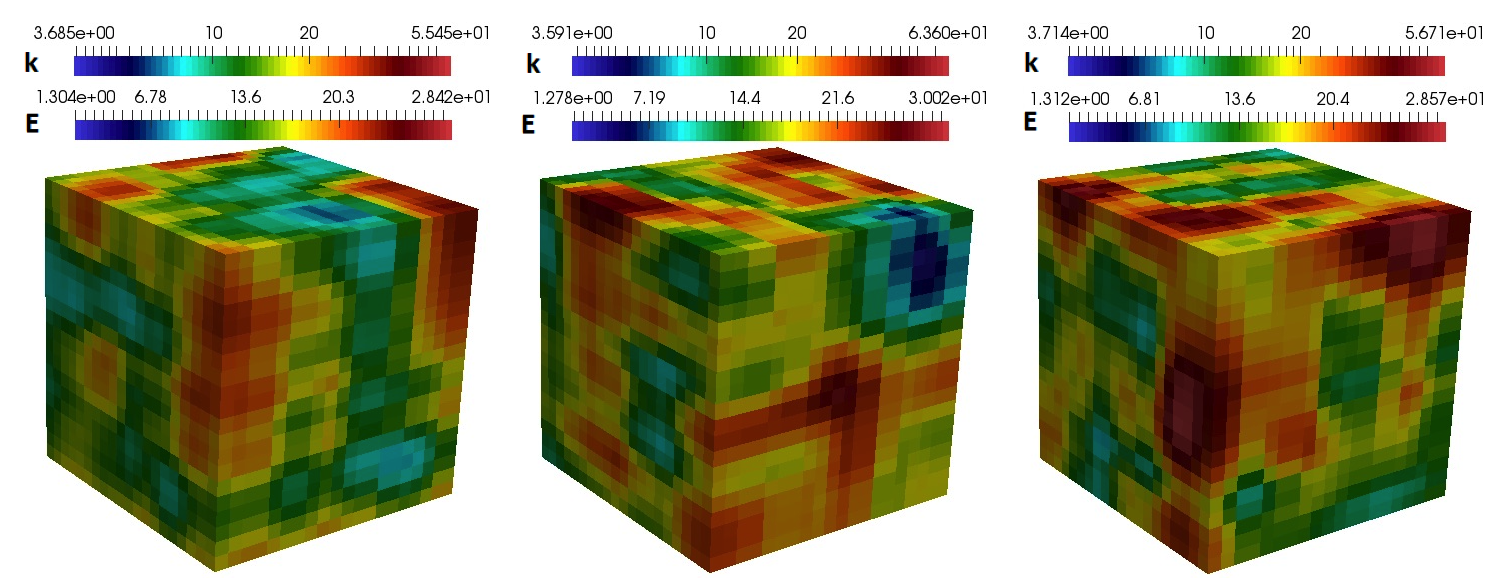}
\caption{Heterogeneous elasticity modulus and permeability. 
First row: Cases 1, 2 and 3 (from left to right) for Test 1. 
Second row: Cases 1, 2 and 3 (from left to right) for Test 2.  
Third row: Cases 1, 2 and 3 (from left to right) for Test 3. }
\label{fig:2d-kE}
\end{figure}

%\begin{figure}[h!]
%\centering
%\includegraphics[width=0.7\linewidth]{homog2d}
%\caption{Distribution of pressure, displacement along $X$ and $Y$ directions at the last moment of time for fine grid (top) and coarse grid using numerical homogenization  (bottom)  for \textit{two-dimensional problem}. Test 1 with case 1}
%\label{fig:2d-u}
%\end{figure}

\begin{table}[h!]
\begin{center}
\begin{tabular}{ | c | c | c | c | c | }
\hline
Case  
& $||e_p||_{1}$  (\%) & $||e_p||_{2}$  (\%) 
& $||e_u||_{1}$ (\%)  & $||e_u||_{2}$  (\%)\\ \hline
\multicolumn{5}{|c|}{Test 1} \\ 
\hline
1 & 2.86139 & 16.3882 & 8.84006 & 16.2718\\
2 & 3.4038 & 15.8993 & 9.21494 & 15.9281\\
3 &  2.82176 & 12.7992 & 8.56077 & 14.3118\\
\hline
\multicolumn{5}{|c|}{Test 2} \\ 
\hline
1 & 3.917 & 11.497 & 3.644 & 12.433 \\
2 & 3.721 & 11.842 & 3.995 & 9.544 \\
3 & 2.939 & 10.849 & 3.712 & 10.216 \\
\hline
\multicolumn{5}{|c|}{Test 3} \\ 
\hline
1 & 7.923 & 20.993 & 8.607 & 24.203 \\
2 & 7.385 & 21.724 & 7.475 & 19.563 \\
3 & 7.223 & 21.603 & 5.399 & 21.777 \\
\hline
\end{tabular}
\end{center}
\caption{Relative $L_2$ and energy errors for displacement and pressure between coarse grid solution and reference solution. Tests 1, 2 and 3 with three cases}
\label{tab:homog}
\end{table}

\textbf{Two-dimensional problem (Tests 1 and 2). }
We solve poroelasticity problem in $\Omega = [0, 1] \times [0,1]$. 
We set following boundary conditions
\[
\begin{split}
u_1 = 0, \quad  \sigma_2 = 0, \quad &x \text{ on } \Gamma_L \\
\sigma_1 = 0,  \quad u_2 = 0, \quad &x \text{ on } \Gamma_B \\
p = 1, \quad &x \text{ on } \Gamma_T.
\end{split}
\]
where $\Gamma_L$ and $\Gamma_R$ are the left and right boundaries, $\Gamma_B$ and $\Gamma_T$ are the bottom and top boundaries,  $\partial \Omega = \Gamma_L \cup \Gamma_R \cup \Gamma_B \cup \Gamma_T$. 

\textbf{Three-dimensional problem (Test 3). } 
We solve model problem in $\Omega = [0, 1] \times  [0, 1] \times [0,1]$.  
We set following boundary conditions
\[
\begin{split}
u_1 = 0, \quad  \sigma_2 = 0, \quad \sigma_3 = 0, 
\quad &x \text{ on } \Gamma_L \\
\sigma_1 = 0,  \quad u_2 = 0, \quad \sigma_3 = 0,  
\quad  &x \text{ on } \Gamma_B \\
\sigma_1 = 0,  \quad \sigma_2 = 0,  \quad u_3 = 0, 
\quad  &x \text{ on } \Gamma_W \\
p = 1, \quad &x \text{ on } \Gamma_T.
\end{split}
\]
where 
$\Gamma_L$ and $\Gamma_R$ are the left and right boundaries, 
$\Gamma_B$ and $\Gamma_T$ are the bottom and top boundaries,  
$\Gamma_F$ and $\Gamma_W$ are the forward and backward boundaries,  $\partial \Omega = \Gamma_L \cup \Gamma_R \cup \Gamma_B \cup \Gamma_T  \cup \Gamma_F \cup \Gamma_B$. 

We compute relative $L_2$ and energy errors between fine grid (reference) solution ($p_f$, $u_f$) and coarse grid solution ($p$, $u$)
\[
||e_p||^2_{1} = 
\frac{\int_{\Omega} (p_f - p)^2 dx }{\int_{\Omega} p_f dx } , \quad 
||e_u||^2_{1} =
\frac{\int_{\Omega} (u_f - u)^2 dx }{\int_{\Omega} u_f^2 dx },
\]\[
||e_p||^2_{2} = 
\frac{\int_{\Omega} k \nabla(p_f - p) \cdot \nabla(p_f - p) \, dx }{\int_{\Omega} k \nabla(p_f ) \cdot \nabla(p_f) \, dx }, \quad 
||e_u||^2_{2} = 
\frac{\int_{\Omega} \sigma(u_f - u) : \varepsilon(u_f - u) \, dx }{\int_{\Omega} \sigma(u_f ) : \varepsilon(u_f) \, dx }. 
\]

We consider three test cases for different heterogeneity (see Figure \ref{fig:2d-kE}) for Tests 1, 2 and 3.
For 2D problems (Tests 1 and 2), coarse grid is $10 \times 10$ and fine grid is $320 \times 320$.
%In Figure \ref{fig:2d-u}, we present the distribution of pressure and displacement along the X and Y directions at the last moment of time for the coarse and fine grids for Test 1 with case 1.
For 3D problems (Test 3), coarse grid is $5 \times 5 \times 5$ and fine grid is $60 \times 60 \times 60$.
In Table \ref{tab:homog}, we present the relative $L_2$ and energy errors of the coarse grid solver for Tests 1,2 and 3.
We observe good accuracy for both 2D and 3D test problems with different heterogeneity.
Next, we consider a machine learning technique for fast prediction of the effective properties for Tests 1, 2 and 3.

\subsection{Results of the learning process}

For dataset construction, we generate a set of random realizations ${Y(x, \omega_l)}$ ($l = \overline{1,M}$) and construct 
\[
k(x,\omega_l)=\exp\left(Y^{}(x,\omega_l)\right),
\quad 
E(x,\omega_l)=\bar{E}(x)+\alpha Y(x,\omega_l).
\]
Next, we define a uniform $N \times N$, divide $k(x, \omega_l) = \cup_i k^{K_i}(x, \omega_l)$ and $E(x,\omega_l) = \cup_i E^{K_i}(x, \omega_l)$ into local domains $K_i$, where $i = \overline{1,N_c}$ and $N_c = N \cdot N$.

Then for each $k^{K_i}(x, \omega_l)$ and  $E^{K_i}(x, \omega_l)$ using  numerical homogenization technique, we calculate $\{k^{K_i}(x,\omega_l) \rightarrow  k^{*,K_i} (\omega_l)\}$ and  $\{E^{K_i}(x,\omega_l) \rightarrow  C^{*,K_i}(\omega_l)\}$ for $i = \overline{1,N_c}$ and $l = \overline{1,M}$. Finally, we normalize dataset and use it for machine learning algorithm construction (each component of the effective tensor is normalized separately).

\begin{figure}[h!]
\centering
\includegraphics[width=1.0 \textwidth]{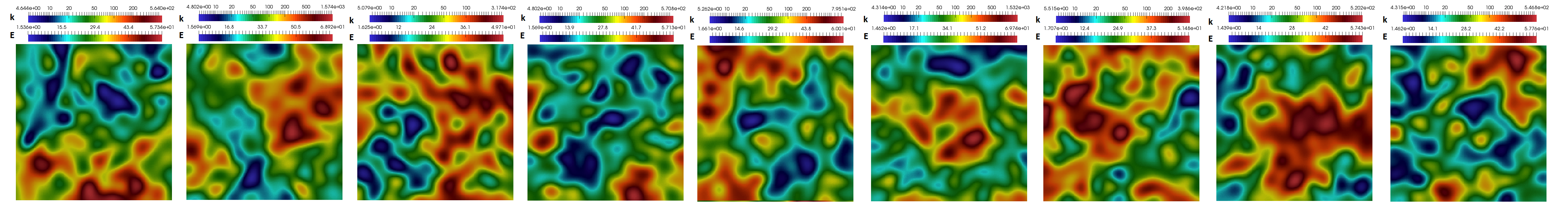}
\caption{Dataset for 2D problem. Test 1}
\label{fig:train2}
\end{figure}

\begin{figure}[h!]
\centering
\includegraphics[width=1.0 \textwidth]{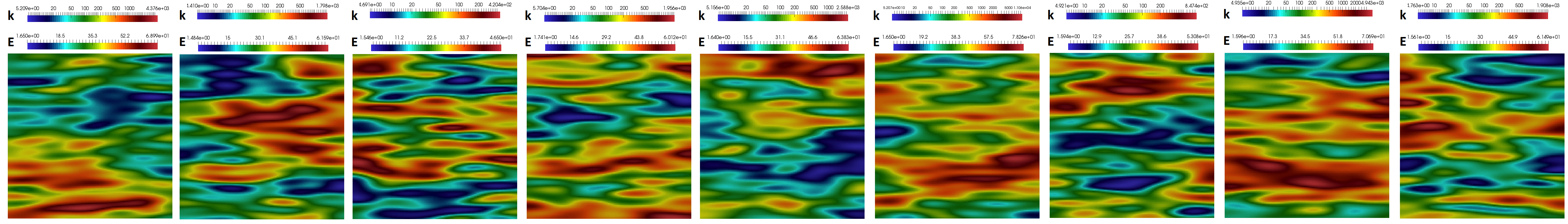}
\caption{Dataset for 2D problem. Test 2}
\label{fig:train2f}
\end{figure}

 \begin{figure}[h!]
\centering
\includegraphics[width=1.0 \textwidth]{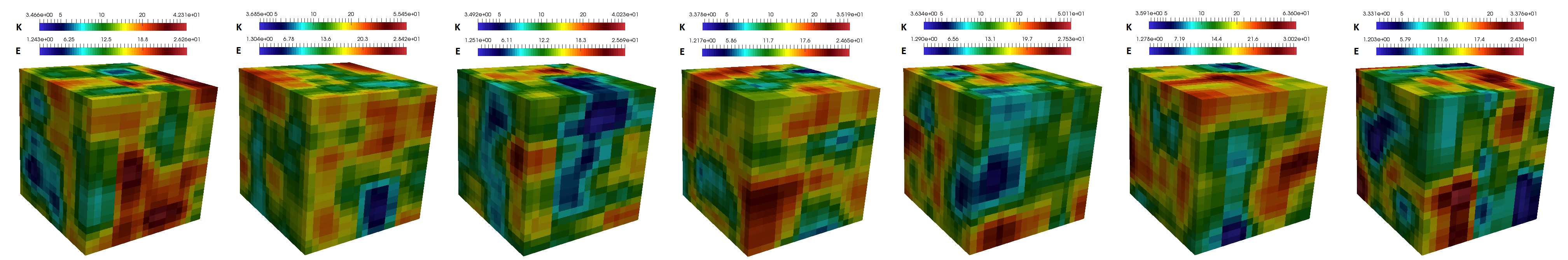}
\caption{Dataset for 3D problem. Test 3}
\label{fig:train3}
\end{figure}

\begin{table}[h!]
\centering
\begin{tabular}{ | c | c | }
\hline  
Input: train and validation set ($i = 1,...,(N_{train} + N_{val})$) & $K_i^{N_l \times N_l}$\\ 
\hline
2D Convolution with $L_1$ filters with $(3,3)$ - kernel and RELU activation &  $L_1 \times K_i^{N_l \times N_l}$  \\ \hline
Max pooling with $(2,2)$ - kernel  &  $L_1 \times K_i^{N_l/2 \times N_l/2}$  \\ 
\hline
2D Convolution with $L_2$ filters with $(3,3)$ - kernel and RELU activation &  $L_2 \times K_i^{N_l/2 \times N_l/2}$  \\ \hline
Max pooling with $(2,2)$ - kernel  &  $L_2 \times K_i^{N_l/4 \times N_l/4}$  \\ 
\hline
2D Convolution with $L_3$ filters with $(3,3)$ - kernel and RELU activation &  $L_3 \times K_i^{N_l/4 \times N_l/4}$  \\ \hline
Max pooling with $(2,2)$ - kernel  &  $L_3 \times K_i^{N_l/8 \times N_l/8}$  \\ 
\hline
2D Convolution with $L_4$ filters with $(3,3)$ - kernel and RELU activation &  $L_4 \times K_i^{N_l/8 \times N_l/8}$  \\ \hline
Max pooling with $(2,2)$ - kernel  &  $L_4 \times K_i^{N_l/16 \times N_l/16}$  \\ 
\hline
Fully connected layer with dropout &  $L_2 \times N_l/16 \times N_l/16$  \\ \hline
Fully connected layer with dropout &  $L_5$  \\ \hline
Output &  $N_{out}$  \\ \hline
\end{tabular}
\caption{Convolutional neural network architecture for prediction effective properties. $N_l = 64$,  $L_1 = 8$, $L_2 = 16$, $L_3 = 32$, $L_4 = 64$, $L_5 = 512$, $N_{out} = 3$ for effective permeability prediction and  $N_{out} = 6$ for elasticity tensor. 2D problem.}
\label{cnn2d}
\end{table}

\begin{table}[h!]
\centering
\begin{tabular}{ | c | c | }
\hline  
Input: train and validation set ($i = 1,...,(N_{train} + N_{val})$) & $K_i^{N_l \times N_l \times N_l}$\\ \hline
3D Convolution with $L_1$ filters with $(3,3,3)$ - kernel and RELU activation &  $L_1 \times K_i^{N_l \times N_l \times N_l}$  \\ \hline
Max pooling with $(2,2,2)$ - kernel  &  $L_1 \times K_i^{N_l/2 \times N_l/2\times N_l/2}$  \\ \hline
3D Convolution with $L_2$ filters with $(3,3,3)$ - kernel and RELU activation &  $L_2 \times K_i^{N_l/2 \times N_l/2 \times N_l/2}$  \\ \hline
Max pooling with $(2,2,2)$ - kernel  &  $L_2 \times K_i^{N_l/4 \times N_l/4 \times N_l/4}$  \\ \hline
Fully connected layer with dropout &  $L_2 \times N_l/4 \times N_l/4 \times N_l/4$  \\ \hline
Fully connected layer with dropout &  $L_3$  \\ \hline
Output &  $N_{out}$  \\ \hline
\end{tabular}
\caption{Convolutional neural network architecture for prediction effective properties. $N_l = 12$,  $L_1 = 16$, $L_2 = 32$, $L_3 = 512$, $N_{out} = 6$ for effective permeability prediction and  $N_{out} = 21$ for elasticity tensor. 3D problem.}
\label{cnn3d}
\end{table}

% TEST 1
\begin{figure}[h!]
\centering
\includegraphics[width=0.32 \textwidth]{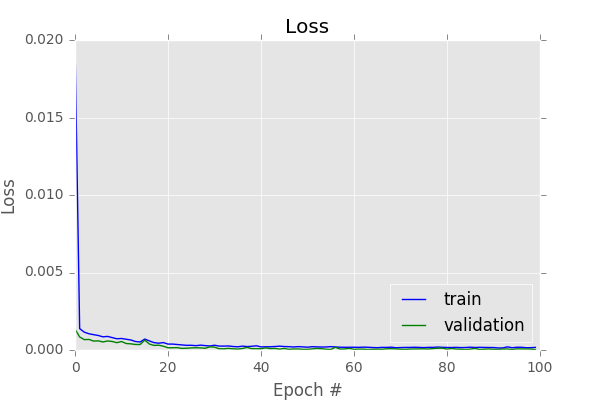}
\includegraphics[width=0.32 \textwidth]{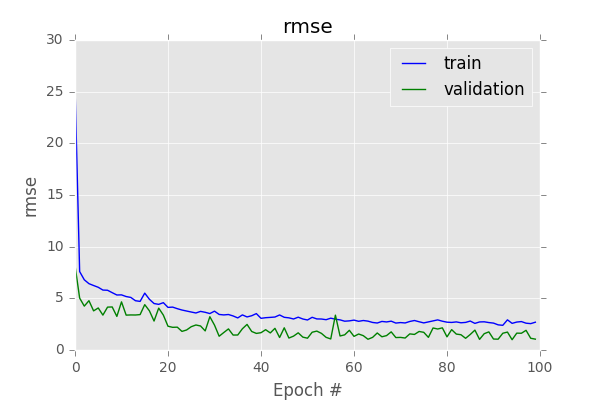}
\includegraphics[width=0.32 \textwidth]{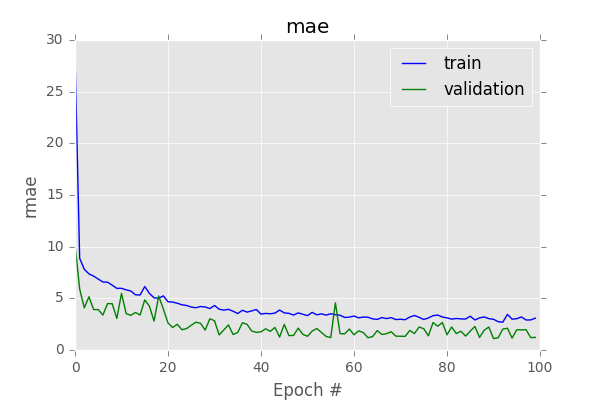}
\caption{Learning process for Test 1. 
Loss, relative root mean square error and relative mean absolute error  vs epoch number for numerical homogenization. 
Effective permeability tensor for Test 1 }
\label{fig:loss1}
\end{figure}

\begin{figure}[h!]
\centering
\includegraphics[width=0.5 \textwidth]{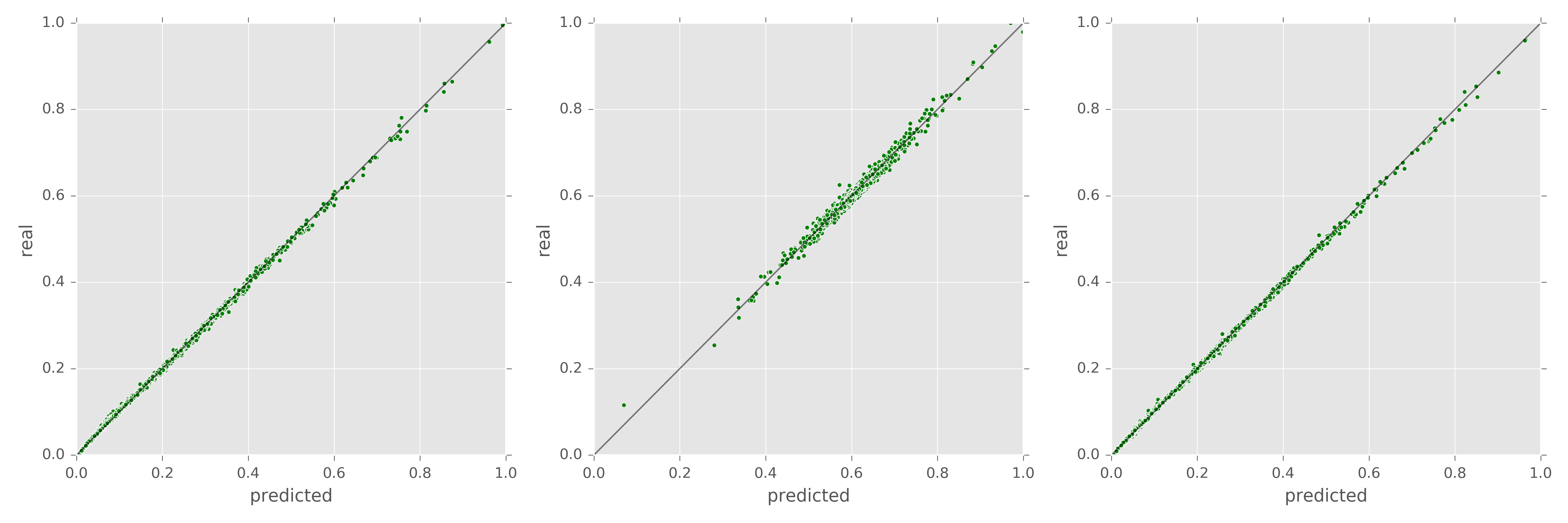}\\
\includegraphics[width=0.5 \textwidth]{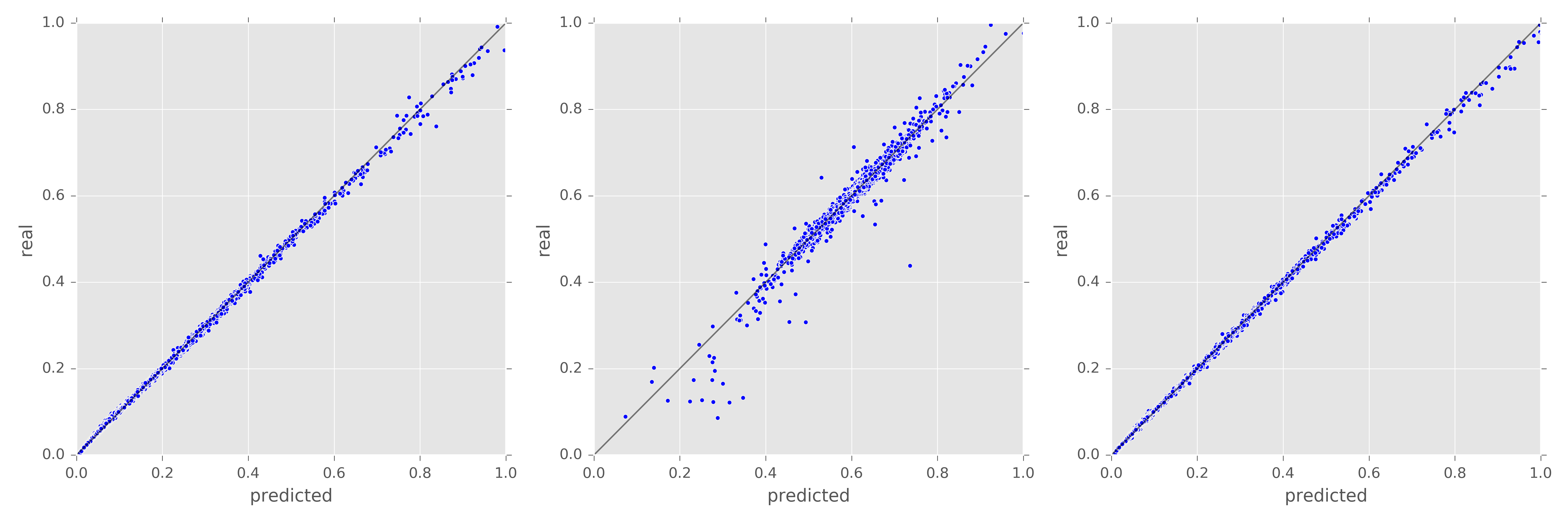}
\caption{Learning performance of CNN for 2D problem (Test 1). 
Effective permeability tensor, $Y_i = \{k^{\star, K_i}_{11},  k^{\star, K_i}_{12}, k^{\star, K_i}_{22} \}$ (from left to right). 
Parity plots comparing preference property values against
predictions made using CNN. 
First row: train  and validation dataset (green color). 
Second row: test dataset (blue color)}
\label{fig:corr2dk}
\end{figure}

\begin{figure}[h!]
\centering
\includegraphics[width=1 \textwidth]{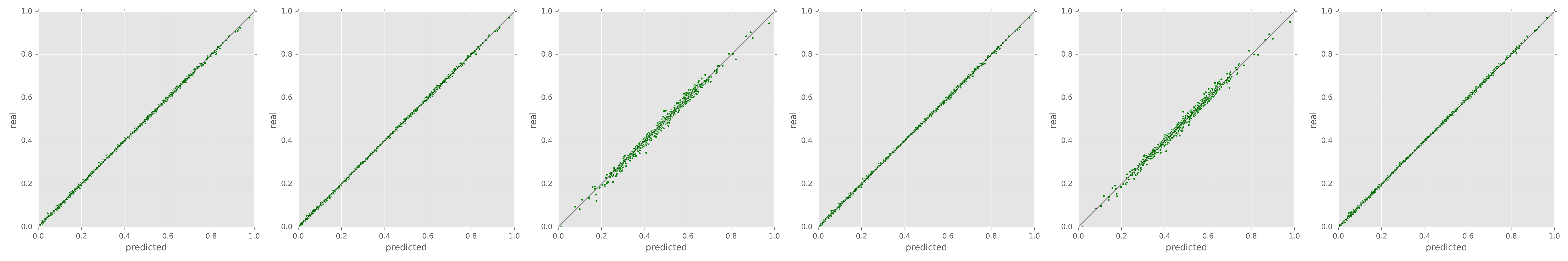}\\
\includegraphics[width=1 \textwidth]{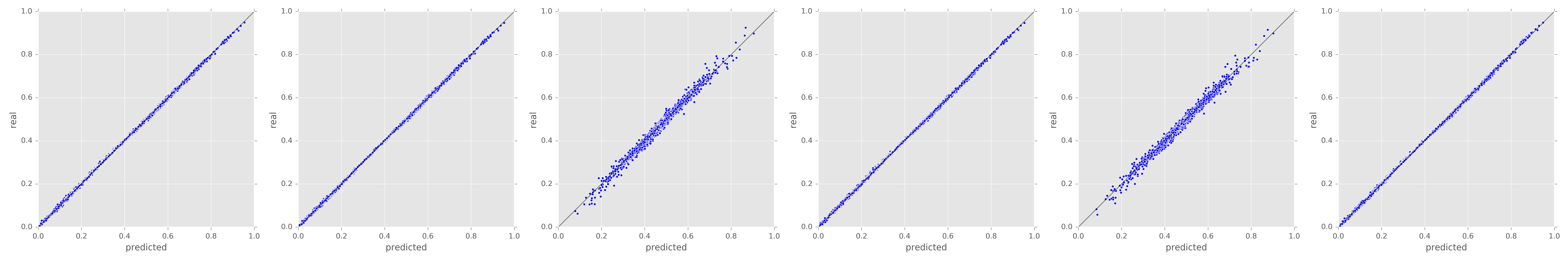}
\caption{Learning performance of CNN for 2D problem (Test 1). 
Effective elasticity tensor, 
$Y_i = \{
C^{\star, K_i}_{1111},  C^{\star, K_i}_{1122}, C^{\star, K_i}_{1112}, C^{\star, K_i}_{2222},  C^{\star, K_i}_{2212}, C^{\star, K_i}_{1212} 
\}$ (from left to right). 
Parity plots comparing preference property values against
predictions made using CNN. 
First row: train and validation dataset (green color). 
Second row: test dataset (blue color)}
\label{fig:corr2dE}
\end{figure}

\begin{table}[h!]
\begin{center}
\begin{tabular}{ | c | c c c | c c c | }
\hline
\multirow{2}{*}{Error} & \multicolumn{3}{|c|}{Train set} & 
\multicolumn{3}{|c|}{Test set} \\ 
& MSE & RMSE  (\%)  & MAE (\%) & MSE & RMSE  (\%)  & MAE (\%) \\ \hline
\multicolumn{7}{|c|}{effective permeability tensor} \\ \hline
$k_{11}$ 	& 0.047 & 2.190 & 2.436 & 0.059 & 2.446 & 2.413 \\
$k_{12}$ 	& 0.012 & 1.131 & 0.834 & 0.074 & 2.730 & 1.083 \\
$k_{22}$ 	& 0.040 & 2.013 & 2.248 & 0.051 & 2.279 & 2.254 \\
\hline
$k$ 				& 0.016 & 1.290 & 1.218 & 0.071 & 2.671 & 1.436 \\ 
\hline
\multicolumn{7}{|c|}{effective elasticity tensor} \\ \hline
$C_{1111}$ 	& 0.008 & 0.903 & 0.788 & 0.008 & 0.942 & 0.819 \\
$C_{1122}$ 	& 0.007 & 0.853 & 0.756 & 0.007 & 0.893 & 0.795 \\
$C_{1112}$ 	& 0.026 & 1.625 & 1.177 & 0.033 & 1.840 & 1.292 \\
$C_{2222}$ 	& 0.008 & 0.931 & 0.832 & 0.009 & 0.965 & 0.867 \\
$C_{2212}$ 	& 0.025 & 1.602 & 1.175 & 0.034 & 1.848 & 1.303 \\
$C_{1212}$ 	& 0.006 & 0.810 & 0.695 & 0.007 & 0.853 & 0.728 \\
\hline
$C$ 					& 0.015 & 1.242 & 0.929 & 0.019 & 1.287 & 0.998 \\
\hline
\end{tabular}
\end{center}
\caption{Learning performance of CNN for 2D problem. Errors for Test 1}
\label{tab:test2d}
\end{table}

% TEST 2
\begin{table}[h!]
\begin{center}
\begin{tabular}{ | c | c c c | c c c | }
\hline
\multirow{2}{*}{Error} & \multicolumn{3}{|c|}{Train set} & 
\multicolumn{3}{|c|}{Test set} \\ 
& MSE & RMSE  (\%)  & MAE (\%) & MSE & RMSE  (\%)  & MAE (\%) \\ \hline
\multicolumn{7}{|c|}{effective permeability tensor} \\ \hline
$k_{11}$ 	& 0.101 & 3.191 & 2.791 & 0.112 & 3.350 & 2.744 \\
$k_{12}$ 	& 0.091 & 3.019 & 1.516 & 0.246 & 4.960 & 1.962 \\
$k_{22}$ 	& 0.168 & 4.110 & 3.634 & 0.290 & 5.387 & 3.701 \\ \hline
$k$ 				& 0.094 & 3.081 & 1.874 & 0.238 & 4.885 & 2.282 \\
\hline
\multicolumn{7}{|c|}{effective elasticity tensor} \\ \hline
$C_{1111}$ 	& 0.010 & 1.014 & 0.836 & 0.010 & 1.015 & 0.815 \\
$C_{1122}$ 	& 0.010 & 1.002 & 0.845 & 0.009 & 0.988 & 0.819 \\
$C_{1112}$ 	& 0.075 & 2.745 & 2.048 & 0.092 & 3.039 & 2.128 \\
$C_{2222}$ 	& 0.009 & 0.988 & 0.806 & 0.009 & 0.988 & 0.796 \\
$C_{2212}$ 	& 0.072 & 2.682 & 1.922 & 0.085 & 2.932 & 1.973 \\
$C_{1212}$ 	& 0.009 & 0.971 & 0.804 & 0.009 & 0.978 & 0.788 \\ \hline
$C$ 					& 0.035 & 1.888 & 1.268 & 0.039 & 1.987 & 1.261 \\
\hline
\end{tabular}
\end{center}
\caption{Learning performance of CNN for 2D problem. Errors for Test 2}
\label{tab:test2df}
\end{table}

% TEST 3
\begin{table}[h!]
\begin{center}
\begin{tabular}{ | c | c c c | c c c | }
\hline
\multirow{2}{*}{Error} & \multicolumn{3}{|c|}{Train set} & 
\multicolumn{3}{|c|}{Test set} \\ 
& MSE & RMSE  (\%)  & MAE (\%) & MSE & RMSE  (\%)  & MAE (\%) \\ \hline
\multicolumn{7}{|c|}{effective permeability tensor} \\ \hline
$k_{11}$ 	& 0.027 & 1.644 & 1.191 & 0.041 & 2.039 & 1.287 \\
$k_{22}$ 	& 0.021 & 1.468 & 1.108 & 0.049 & 2.235 & 1.263 \\
$k_{33}$ 	& 0.025 & 1.586 & 1.181 & 0.051 & 2.274 & 1.329 \\
$k_{12}$ 	& 0.040 & 2.003 & 1.262 & 0.118 & 3.446 & 1.603 \\
$k_{13}$ 	& 0.050 & 2.238 & 1.351 & 131 & 3.629 & 1.642 \\
$k_{23}$ 	& 0.036 & 1.904 & 1.372 & 0.207 & 4.559 & 1.800 \\ \hline
$k$ 				& 0.039 & 1.991 & 1.283 & 0.132 & 3.646 & 1.567 \\
\hline
\multicolumn{7}{|c|}{effective elasticity tensor} \\ \hline
$C_{1111}$ 	& 0.017 & 1.337 & 1.063  & 0.020   & 1.434  & 1.096  \\ 
$C_{1122}$ 	& 0.016 & 1.269 & 1.015  & 0.018   & 1.343  & 1.042  \\ 
$C_{1133}$ 	& 0.016 & 1.272 & 1.004  & 0.019   & 1.381  & 1.047   \\
$C_{1112}$ 	& 0.055 & 2.361 & 1.675  & 0.079   & 2.823  & 1.926   \\
$C_{1123}$   & 0.028 & 1.702 & 1.172  & 0.052   & 2.299  & 1.435  \\
$C_{1131}$ 	& 0.039 & 1.985 & 1.420  & 0.059   & 2.447  & 1.630   \\
$C_{2222}$ 	& 0.014 & 1.189 & 0.923  & 0.016   & 1.282  & 0.960   \\
$C_{2233}$ 	& 0.014 & 1.216 & 0.973  & 0.017   & 1.307  & 1.007   \\
$C_{2212}$ 	& 0.049 & 2.220 & 1.559  & 0.074   & 2.735  & 1.811   \\
$C_{2223}$ 	& 0.031 & 1.781 & 1.251  & 0.056   & 2.378  & 1.524   \\
$C_{2231}$   & 0.038 & 1.950 & 1.373  & 0.058   & 2.418  & 1.584  \\
$C_{3333}$   & 0.015 & 1.262 & 0.981  & 0.018   & 1.375  & 1.028  \\
$C_{3312}$ 	& 0.051 & 2.261 & 1.593  & 0.079   & 2.813  & 1.854   \\
$C_{3323}$ 	& 0.027 & 1.664 & 1.149  & 0.050   & 2.256  & 1.419   \\
$C_{3331}$ 	& 0.038 & 1.966 & 1.394  & 0.059   & 2.444  & 1.603   \\
$C_{1212}$ 	& 0.017 & 1.309 & 1.013  & 0.020   & 1.418  & 1.055   \\
$C_{1223}$ 	& 0.059 & 2.429 & 1.717  & 0.085   & 2.929  & 1.933   \\
$C_{1231}$   & 0.042 & 2.070 & 1.469  & 0.073   & 2.716  & 1.780  \\
$C_{2323}$ 	& 0.014 & 1.218 & 0.958  & 0.017   & 1.318  & 0.989   \\
$C_{2331}$ 	& 0.081 & 2.853 & 2.009  & 0.115   & 3.396  & 2.275   \\
$C_{3131}$ 	& 0.015 & 1.229 & 0.995  & 0.017   & 1.331  & 1.027   \\ \hline
$C$ 					& 0.033 & 1.838 & 1.286  & 0.050   & 2.255  & 1.452  \\
\hline
\end{tabular}
\end{center}
\caption{Learning performance of CNN for 3D problem. Errors for Test 3}
\label{tab:test3d}
\end{table}

We set $M = 100$, $N = 10$ for 2D problems and $N = 5$ for 3D problem. Therefore, we have dataset with size $L = 100 \cdot 10^2$ (for Tests 1 and for Test 2, separately) and $L = 100 \cdot 5^3$ (Test 3). 
Therefore input of the network is  $d$ - dimensional normalized array 
\[
X_j = k^{K_i}(x,\omega_l) \quad  \text{ or } \quad 
X_j = E^{K_i}(x,\omega_l), \quad 
j = l \cdot M + i, 
\]
where $X_j$ has size $N_l^d$, where $N_l = 64$ (Tests 1 and 2) and $N_l =12$ (Test 3).  

The output for the numerical homogenization is the  normalized array of the effective properties
\[
Y_j = k^{*,K_i}(x,\omega_l) \quad  \text{ or } \quad 
Y_j = C^{*,K_i} (x,\omega_l), \quad 
j = l \cdot M + i, \quad 
i = \overline{1,N_c}, \quad 
l = \overline{1,M}, 
\]
where $k^{*,K_i}$ and $C^{*,K_i}$ are tensors that represented as array with symmetric property. 

The dataset  $\{X_j \rightarrow Y_j\}$ used to train, validate and test the neural network (see Figures \ref{fig:train2}, \ref{fig:train2f}  and \ref{fig:train3} for illustration of $k(x,\omega_l)$ and $E(x,\omega_l)$). 
Each dataset is divided into 3200:800:6000  ratio for training, validation and test sets.   

% ml
For calculations, we use 100 epochs with a batch size $N_c$ (number of coarse cells) and Adam optimizer with  learning rate $\epsilon = 0.001$. 
For accelerating of the training process of the CNN, we use GPU (Nvidia GeForce GTX 1080 Ti). 
In order to prevent overfitting, we use dropout with rate 10 \%.  
We use $d$ - dimensional convolutions and max pooling layers with size $3^d$ and  $2^d$, respectively. The architectures of the CNN for 2D and 3D problems are presented in Tables \ref{cnn2d} and \ref{cnn3d}. 
Proposed CNN for 2D problems contains 11 layers and 7 layers for 3D problems. The input data $X_i$ is a $N_l^d$ array with $d = 2,3$, $N_l = 64$ for 2D and $N_l = 12$ for 3D. %Convolution layers with RELU activation have size $3^d$ and alternate with pooling layers with size $2^d$. 
Convolution layer contains $L_i$ feature maps.  Finally, we apply three fully connected layers.

For error calculation on the train and test dataset, we use  mean square errors, relative mean absolute and relative root mean square errors 
\[
MSE = \sum_i |Y_i - \tilde{Y}_i|^2,
\quad  
MAE = \frac{\sum_i |Y_i - \tilde{Y}_i|  }{\sum_i |Y_i|},
\quad 
RMSE = \sqrt{ \frac{\sum_i |Y_i - \tilde{Y}_i|^2 }{\sum_i |Y_i|^2 } },
\]
where $Y_i$  and $\tilde{Y}_i$ denotes the reference and predicted values for sample $X_i$

Convergence of the loss function for Test 1 (2D)  is presented in Figure \ref{fig:loss1}, where we plot the relative root mean square error (RMSE) and relative mean absolute error (MAE) vs epoch number. We depict loss functions for train and validation sets for learning of the effective permeability tensor. 
In Figures \ref{fig:corr2dk} and \ref{fig:corr2dE}, we present  a parity plots comparing preference property values against predicted using CNN. %Convergence of the loss function, is depicted in Figure \ref{fig:loss1}. 
Learning performance of CNN for 2D and 3D problems are presented in Tables \ref{tab:test2d}, \ref{tab:test2df} and \ref{tab:test3d} for Tests 1, 2 and 3, respectively. We observe good convergence for the relative errors  for train and test sets.

\subsection{Numerical homogenization with machine learning approach} 

Next, we consider errors between solution of the coarse grid problem with reference  and predicted effective properties. We present results for the Tests 1, 2 and 3. 

\begin{figure}[h!]
\centering
Test 1. \\
\includegraphics[width=0.7 \textwidth]{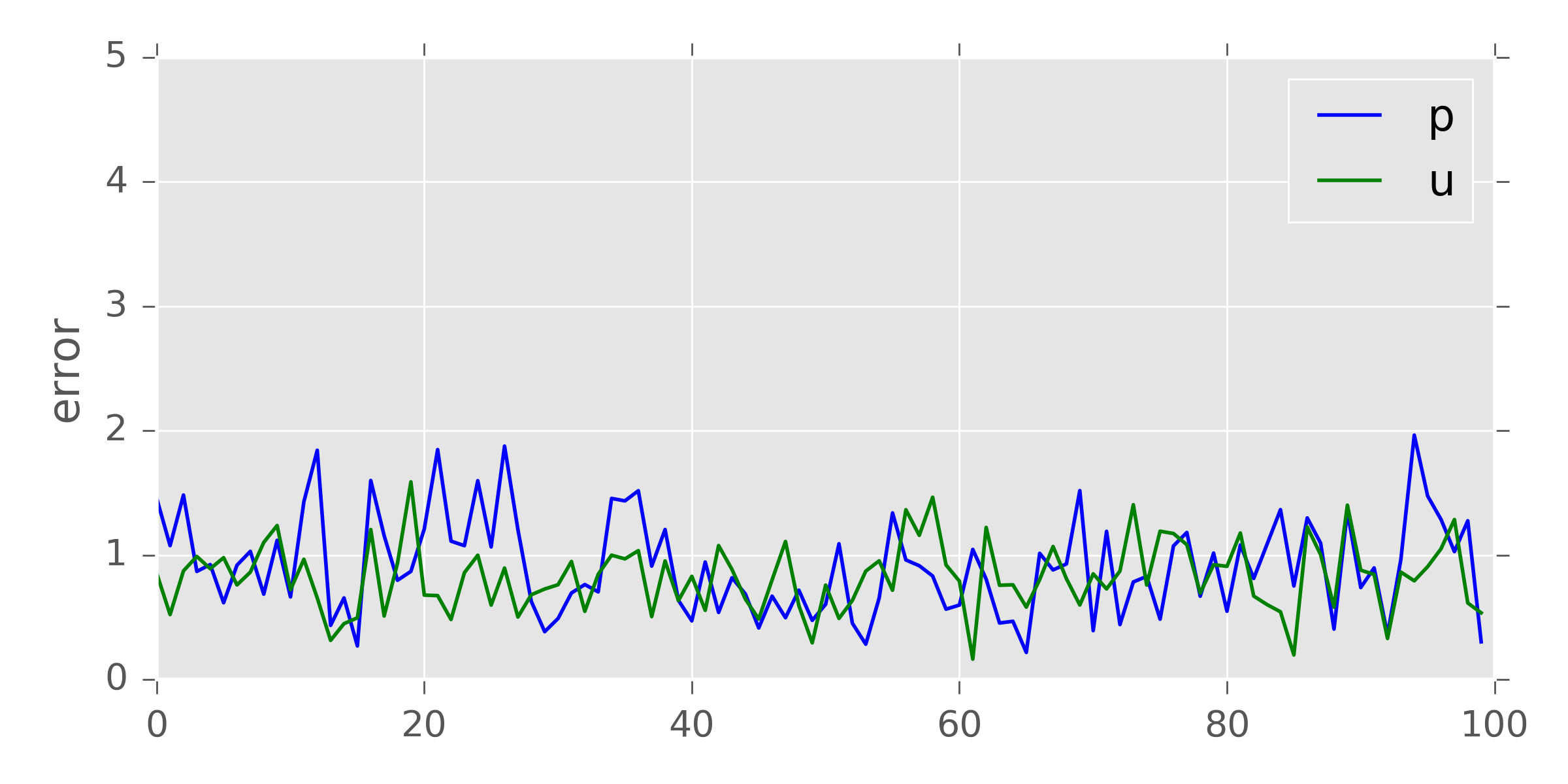}\\
Test 2. \\
\includegraphics[width=0.7 \textwidth]{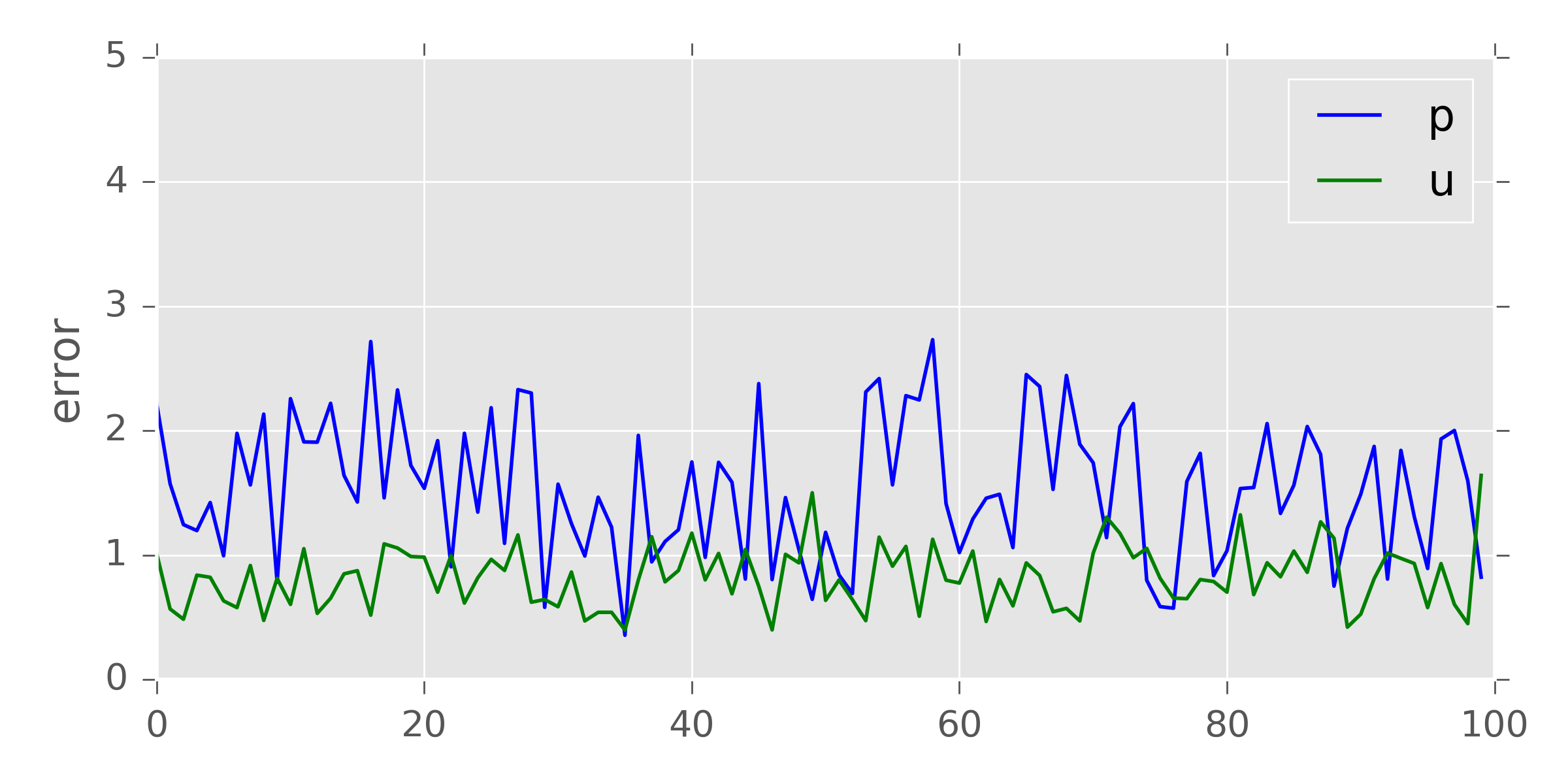}\\
Test 3. \\
\includegraphics[width=0.7 \textwidth]{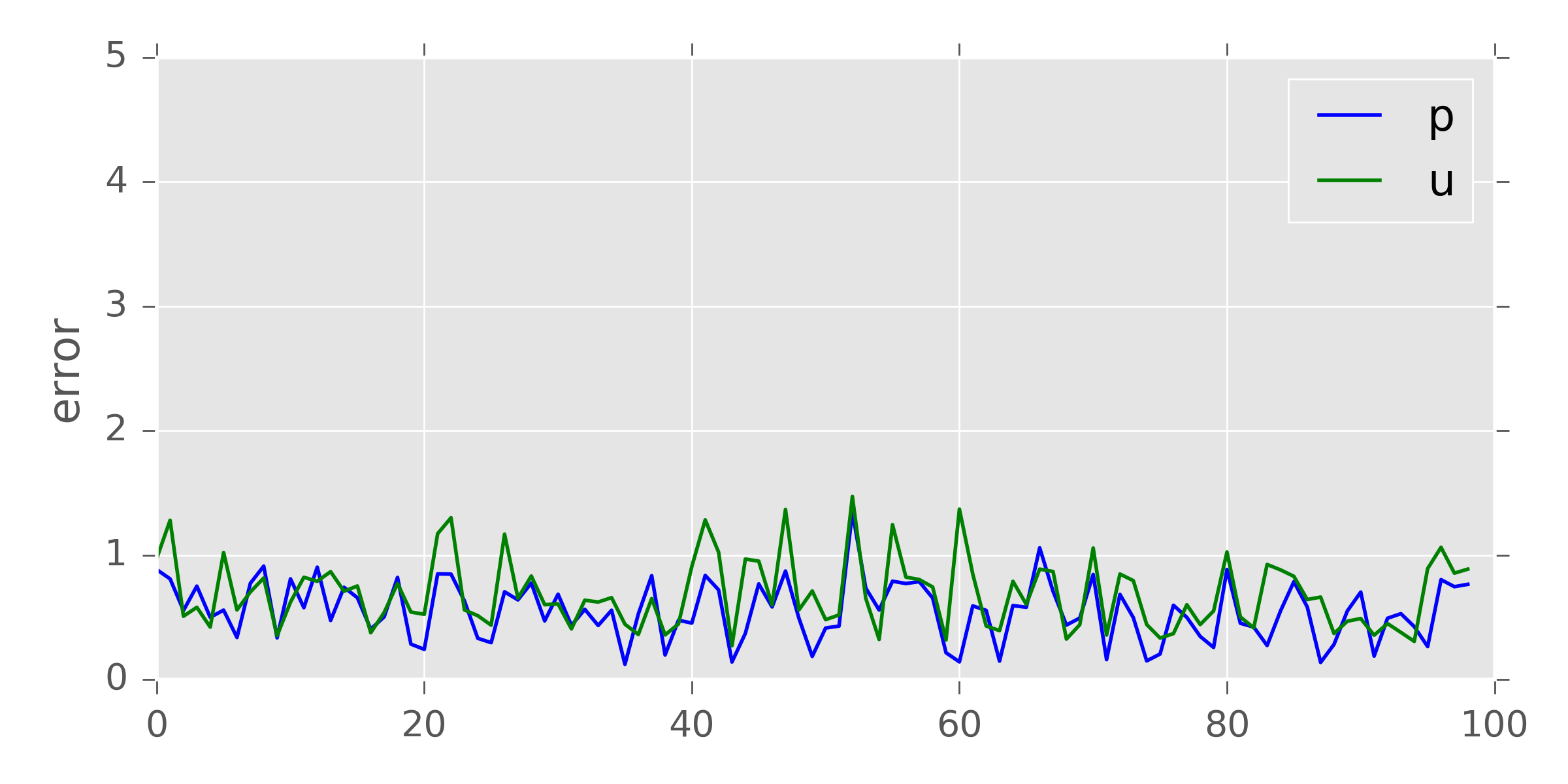}
\caption{Numerical homogenization error with effective coefficients predicted using machine learning algorithm. Test 1, 2 and 3 (from top to bottom). Pressure (blue) and displacement (green) }
\label{fig:err-ml}
\end{figure}

\begin{table}[h!]
\begin{center}
\begin{tabular}{ | c | c | c c | c | c | }
\hline
Test & Offline (GPU) & Online (loading)  & Online (prediction) 
&  Direct solve &  Speedup \\ \hline
\multicolumn{6}{|c|}{effective permeability tensor} \\ 
\hline
Test 1 (2D) 	& 40.688 & 0.547 & 0.171 & 13.025 & $\times$ 76.1  \\
\hline
Test 2 (2D) 	& 40.798 & 0.580 & 0.171 & 13.025 & $\times$ 76.1 \\
\hline
Test 3 (3D) 	& 41.163 & 0.431 & 0.552 & 44.833 & $\times$ 81.2  \\
\hline
\multicolumn{6}{|c|}{effective elasticity tensor} \\ \hline
Test 1 (2D) 	& 42.807 & 0.485 & 0.171 & 23.333 & $\times$ 136.4 \\
\hline
Test 2 (2D) 	& 47.169 & 0.433 & 0.173 & 23.333 & $\times$ 134.8 \\
\hline
Test 3 (3D) 	& 41.789 & 0.484 & 0.532 & 153.82 & $\times$ 289.1 \\
\hline
\end{tabular}
\end{center}
\caption{Time in seconds. Offline time is the time of CNN training and validation on GPU. Online calculations contains neural network loading and prediction of the coarse grid effective properties for one $\Omega$. Direct solve is the time of effective properties calculation for one $\Omega$. Speedup is  Direct solve/Online (prediction) }
\label{tab:test2df}
\end{table}

In Figure \ref{fig:err-ml}, we present results for numerical homogenization for 100 samples of the random field. We depict a relative mean square error in percentages for pressure and displacements with effective coefficients predicted using machine learning algorithm for Test 1, 2 and 3. We observe small errors (near one percent) with fast calculations using a machine learning algorithm. 
%Note that, for porous media with high-construct properties this methodology cannot provide good results because numerical homogenization is inaccurate and multiscale multicontinuum methods should be used. Multiscale multicontinuum methods identify multiple important modes, accurately separate variables and provide good accuracy for coarse grid approximations. We will consider it in future works.

Finally, we discuss the computational gain of the machine learning method that achieved by the fast calculations of the effective properties and local matrices.
We divide calculation on the offline and online stages for machine learning algorithm. On the online stage, we train neural network on the GPU by a given train and validation datasets.
We note that, here we didn't consider time of the dataset construction.
On the offline stage, we have two steps: loading of the preconstructed neural network and prediction of the effective properties for a given fine scale distribution. We compare time of the prediction vs time of direct solution of the local problems for effective properties calculations for a some givens fine scale properties in domain $\Omega$.

We observe high speedup of the calculations using GPU.  Time is 1574.910 seconds if we train the neural network on CPU (2.9 GHz Intel Core i5) and 13.025 seconds on GPU for Test 1 (effective permeability tensor). Therefore, we obtain $\times 38$ faster neural network construction on the GPU (Nvidia GeForce GTX 1080Ti). 
In Table \ref{tab:test2df}, we present time of offline and online stages.
On the table, we depict the offline calculation speedup using direct simulations and prediction using CNN. We have approximately $\times 80$ faster calculations for effective permeability and $\times (130-290)$ faster for effective elasticity tensor.
 
\section{Conclusion} 

In this work, we  considered a numerical homogenization of the poroelasticity problem with stochastic properties. For accelerating of the calculations of the effective properties, we constructed a machine learning algorithm. 
We constructed a machine learning algorithm through convolutional neural network (CNN) to learn a map between input stochastic fields to effective properties. 
We trained neural network on the set of the selected realizations of the local microscale stochastic fields and macroscale characteristics (permeability and elasticity tensors). 
Proposed method is used to make fast and accurate effective property predictions for numerical solution of the poroelasticity problems in stochastic media.

%\section{Acknowledgement}
%Work is supported by the mega-grant of the Russian Federation Government (N 14.Y26.31.0013). 

\bibliographystyle{plain}
\bibliography{lit}

\end{document}